\documentclass[11pt,twoside]{article}
\usepackage{amssymb,amsfonts,amsmath,amsthm,cite,color}
\usepackage{dsfont}
\usepackage{epsfig}
\usepackage{indentfirst}
\usepackage{url}
\usepackage{xspace}
\usepackage{latexsym}
\usepackage{epsfig}
\usepackage{fancybox}
\usepackage{amssymb}
\usepackage{amsmath}
\usepackage{amsfonts}
\usepackage{rotating}
\usepackage{booktabs}
\usepackage{amsthm,stmaryrd}
\usepackage{indentfirst}
\usepackage{color}
\usepackage{geometry}
\usepackage{tikz}
\usepackage{graphicx}
\usepackage{wrapfig}
\usepackage{fancyhdr}
\usepackage{fullpage}
\usepackage{hyperref}
\newtheorem{thm}{Theorem}[section]

\newtheorem{lem}[thm]{Lemma}

\newtheorem{defn}[thm]{Definition}
\newtheorem{rem}[thm]{Remark}
\newtheorem{exam}[thm]{Example}

\makeatletter \@addtoreset{equation}{section}
\def\mid{{\,|\,}}

\def\pf{\noindent {\it Proof.\ }}
\def\qed{\hfill \rule{4pt}{7pt}}

\def\bZ{\mathbb Z}
\def\bN{\mathbb N}

\def\la{\langle}
\def\ra{\rangle}
\newcommand{\bK} { {\mathbb{K}}}

\title{An Algorithmic Approach to the $q$-Summability Problem of Bivariate Rational Functions
}
\author{Rong-Hua Wang
  }
\date{\small
  School of Mathematical Sciences, Tiangong University,
  Tianjin, 300387, China\\
  {\sf wangronghua@tjpu.edu.cn}
}

\begin{document}
\maketitle
%\begin{center}
%
%
%Qing-Hu Hou$^1$ and Rong-Hua Wang$^2$\\[7pt]
%
%$^{1}$School of Mathematics \\
%Tianjin University, Tianjin 300072, P. R. China \\[5pt]
%
%$^{2}$School of Mathematics \\
%Tianjin Polytechnic University, Tianjin 300387, P. R. China
%
%\vskip 3mm
%
% E-mail:
% $^1$qh\_hou@tju.edu.cn ,
% $^2$wangronghua@tjpu.edu.cn
%\end{center}
\begin{abstract}
In 2014, Chen and Singer solved the summability problem of bivariate rational functions.
Later an algorithmic proof was presented by Hou and the author.
In this paper, the algorithm will be simplified and adapted to the $q$-case.
\end{abstract}

\section{Introduction}
Symbolic summation is a classical topic in combinatorics and mathematical physics.
One of the central problems in symbolic summation is to decide whether a given sum can be expressed in a closed form, which was fully answered by Gosper's algorithm~\cite{Gosper1978} for indefinite summations of hypergeometric terms.
Based on Gosper's algorithm, Zeilberger~\cite{Zeilberger1990c,Zeilberger1991}
designed a new algorithm to find recurrence relations for single sums of hypergeometric terms, which is known as Zeilberger's algorithm or the method of creative telescoping.
Gosper's and Zeilberger's algorithms occupy a central position in the study of mechanical proofs of combinatorial identities.
The crucial step of both algorithms is to decide whether a given term $T(n)$ can be written as the difference of another term.
If such a term exists, $T(n)$ is said to be summable.
Deciding whether a given term is summable or not is the so-called \emph{summability problem}.

For univariate functions, the summability problem
has been solved rather successfully.
For example, Abramov~\cite{Abramov1975,Abramov1995b} solved the summability problem for rational functions.
Gosper's algorithm~\cite{Gosper1978} settles the summability problem for hypergeometric terms and was later generalized to the D-finite case by Abramov and van Hoeij\cite{abramov_vanhoeij99} and to the difference-field setting by Karr\cite{Karr1981,Karr1985}.

Passing from the univariate case to the multivariate case, the summability problem becomes much more complicated.
Significant progress has been made by Apagodu and Zeilberger
 \cite{Apagodu2006}, Koutschan \cite{Koutschan2010},
Schneider~\cite{Schneider2005} and
Chen et al.~\cite{ChenHouMu2006}.
However, they did not provide a complete answer to
the summability problem of bivariate functions.
The first necessary and sufficient condition for the summability of bivariate functions was presented by Chen and Singer~\cite{ChenSinger2014} for the rational case, and later extended to the remaining mixed cases by Chen in \cite{chen2018ACA}.
Based on the theoretical criterion given in \cite{ChenSinger2014}, Hou and the author \cite{HouWang2015} presented a new criterion and an algorithm for deciding the summability of bivariate rational functions.
The bivariate summability criteria and their mixed analogues
are crucial for solving the existence problems  of telescopers for rational functions in three variables \cite{CHLW2016ISSAC,CDZ2019ISSAC}.

In this paper, the algorithm for detecting the summability of bivariate rational functions will be adapted to the $q$-case.
To a large extent, the $q$-case is analogous to the ordinary case.
The main idea in both situations is to decompose a given function according to different orbits by partial fraction decompositions.
However, in order to obtain concise criteria on the $q$-summability, we have to modify the definition of orbits and thus the whole discussion process and results.
We also provide a much easier proof of the main theorem which provides a criterion on the $q$-summability of rational fractions.
Besides, we show that when a bivariate rational function $f$ is $(\tau_x,\tau_y)$-summable, then there always exist $g,h$ in \emph{reduced forms} such that $f=\tau_xg-g+\tau_yh-h$, where $\tau_x$ and $\tau_y$ denote $q$-shift operators in $x$ and $y$, respectively.

For the sake of readability, we recall some notations and definitions which will be
used frequently.
Throughout the paper we let $\mathbb{K}$ be a field of characteristic zero, $\mathbb{K}(x,y)$ be the field of rational functions in $x,y$ over $\mathbb{K}$ and $q\in\mathbb{K}\setminus\{0\}$ not a root of unity.
Choosing the pure lexicographic order $x\prec y$, a polynomial in $\bK[x,y]$ is called \emph{monic} if its leading coefficient is one.
For a nonzero polynomial $p\in\bK[x,y]$, its degree with respect to the variable
$v\in\{x,y\}$ is denoted by $\deg_v(p)$. We will follow the convention that $\deg_v(0)=-\infty$.

We define $q$-shift operators $\tau_x$ and $\tau_y$ on $\bK(x,y)$ as
\begin{equation}\label{eq:q-shift operator}
\tau_xf(x,y)=f(qx,y)\quad \hbox{and}\quad \tau_yf(x,y)=f(x,qy)
\end{equation}
for all $f\in \bK(x,y)$.
A rational function $f\in \bK(x,y)$ is said to be \emph{$(\tau_x,\tau_y)$-summable}, abbreviated as
\emph{$q$-summable} in certain instances, if there exist $g,h\in \bK(x,y)$ such that
\[f=\tau_xg-g+\tau_yh-h.\]
When $h=0$ or $g=0$, $f$ is called \emph{$\tau_x$-summable} or \emph{$\tau_y$-summable}, respectively.
Deciding whether a given rational function in $\bK(x,y)$ is $q$-summable is the so-called \emph{$q$-summability problem}, which is the central problem to be solved in this paper.

The rest of the paper is organized as follows.
Firstly we extend the concept of bivariate dispersion set to the $q$-case and present a decision procedure to determine whether two multivariate polynomials are $q$-shift equivalent in Section~\ref{q-dispersion}.
Then in Section~\ref{summability} we reduce the $q$-summability problem of a general bivariate rational function to that of a much simpler one, which can be determined by testing the $q$-shift equivalence and solving a $q$-difference equation.
Finally, Section \ref{equation} is devoted to solving the $q$-difference equation.

\section{Testing the $q$-Shift equivalence}\label{q-dispersion}
Let $G=\la\tau_x,\tau_y\ra$ be the free Abelian multiplicative group generated by $\tau_x$ and $\tau_y$.
Let $f\in\mathbb{K}[x,y]$ and $H$ be a subgroup of $G$.
We call
\[
  [f]_H:=\{q^{v}\cdot\sigma(f)|v\in\mathbb{Z},\sigma\in H\}
\]
the \emph{$H$-orbit} at $f$.
Two polynomials $f,g\in\mathbb{K}[x,y]$ are said to be \emph{$H$-equivalent} if
$[f]_H=[g]_H$, denoted by $f\sim_Hg$.
Apparently the relation $\sim_H$ is an equivalence relation.
We will just say $f,g$ are \emph{$q$-shift equivalent} when $f\sim_Hg$ for some $H$ and $H$ is clear from the context.

It should be noticed that, the definitions of $H$-orbit here and the ordinary case as in~\cite{HouWang2015} are different.
The reason that we introduce a power of $q$ is to simplify the $q$-summability problem by decomposing a given rational function according to different orbits by partial fraction decomposition.
Due to this difference, the process for checking the $q$-summability of a rational function is different from that of the summability problem.

In this section, we will present an algorithm to determine whether two polynomials are $q$-shift equivalent, which is also called the \emph{$q$-shift equivalence testing problem.} In fact, this problem can be solved more generally.

Let $\bK[x_1,\ldots,x_n]$ be the ring of polynomials in $x_1,\ldots,x_n$ over $\bK$ and
$q$-shift operators $\tau_{x_i},i=1,2,\ldots,n$ be defined similarly as \eqref{eq:q-shift operator}.
\begin{defn}
For $f,g\in\bK[x_1,\ldots,x_n]$, the $q$-\emph{dispersion set} of $f$ and $g$ is defined as
\[
qDisp(f,g)=\{(\ell,\ell_1,\ldots,\ell_n)\in\mathbb{Z}^{n+1}\ |\ f=q^{\ell}\tau_{x_1}^{\ell_1}\cdots\tau_{x_n}^{\ell_n}(g)\}.
\]
\end{defn}
From the above definition, we know $qDisp(f,g)=\emptyset$ if and only if $f$ and $g$ are not $q$-shift equivalent.
Next we will show that $q$-dispersion set of any two polynomials is computable, which solves the $q$-shift equivalence testing problem.
\begin{thm}\label{th:q-shit equivalent}
Given any two polynomials $f,g \in \bK[x_1, ..., x_n]$, we can determine qDisp(f,g).
\end{thm}
\pf Let $T(f),T(g)$ be the set of all nonzero monic monomials appearing in $f$ and $g$, respectively.
Since $q$-shift operators do not change the term structure, it is easy to see $qDisp(f,g)=\emptyset$ if $T(f)\neq T(g)$.
Suppose $T(f)=T(g)=T$,
\[
f=\sum_{\mathbf{x}^\mathbf{m}\in T}a_{\mathbf{m}}\mathbf{x}^\mathbf{m} \quad
\text {  and  } \quad
g=\sum_{\mathbf{x}^\mathbf{m}\in T}b_{\mathbf{m}}\mathbf{x}^\mathbf{m},
\]
where $\mathbf{x}=(x_1,\dots,x_n),\ \mathbf{m}=(m_1,\ldots,m_n)\in\mathbb{Z}^n$ and
$\mathbf{x}^\mathbf{m}=x_1^{m_1}\cdots x_n^{m_n}$.

Assuming $f=q^{\ell}\tau_{x_1}^{\ell_1}\cdots\tau_{x_n}^{\ell_n}g$ for some $\ell,\ell_i\in\mathbb{Z}$.
For any $\mathbf{x}^\mathbf{m}\in T$, comparing the coefficients both sides leads to
\begin{equation}\label{eq:qDisp set}
a_{\mathbf{m}}=q^{\ell+m_1\ell_1+\cdots+m_n\ell_n}b_{\mathbf{m}},
\end{equation}
which is impossible when $b_{\mathbf{m}}/a_{\mathbf{m}}$ is not a power of $q$ as $q$ is not a root of unity.
Thus we can assume $b_{\mathbf{m}}/a_{\mathbf{m}}=q^m$ for some $m\in\bZ$, then Equation~\eqref{eq:qDisp set} leads to
\begin{equation}\label{eq:Hermite normal form}
\ell+m_1\ell_1+\cdots+m_n\ell_n+m=0,
\end{equation}
which is a linear Diophantine equation in unknowns $\ell$ and $\ell_{i},i=1,...,n$.
The arbitrariness of $\mathbf{x}^\mathbf{m}\in T$ leads to a
linear equation system which can be solved by the computation of Hermite normal form~\cite{kannan1979polynomial} of matrices.\qed

\begin{exam}
Let
\[
f=x+y+1 \quad \hbox{and}\quad g=x+y+q.
\]
It is easy to see $T(f)=T(g)=\{x,y,1\}$.
Suppose $f=q^{\ell}\tau_x^{\ell_1}\tau_y^{\ell2}g$.
Then comparing the coefficients on both sides lead to
\[
\ell+\ell_1=0, \quad \ell+\ell_2=0\quad \text{and}\quad \ell+1=0.
\]
Solving this linear system over $\mathbb{Z}$, we get $\ell=-1$ and $\ell_1=\ell_2=1$.
Hence $qDisp(f,g)=\{(-1,1,1)\}$ and $f = q^{-1}\tau_{x} \tau_{y}(g)$.
\end{exam}
The following example shows $qDisp(f,g)$ can also be infinite.
\begin{exam}
Suppose $n$ is a positive integer and $f=x^n+y^n$.
It is easy to check that $qDisp(f,f)=\{(-nm,m,m)\mid m\in\bZ\}$ is infinite.
\end{exam}
\section{Criteria on $q$-summability}\label{summability}

This section is devoted to determining whether a given bivariate rational function is $q$-summable.
Since we can reduce the $q$-summablity problem of bivariate rational functions to that of univariate ones, we first demonstrate how to determine whether a  given rational function $f$ is $\tau_y$-summable using the $q$-polynomial residue.

\subsection{Univariate case}\label{SE:UniqSumCri}

%For $f,g\in\bK[x]$, we define their $q$-\emph{dispersion set}
%to be
%\[
%qDisp(f,g)=\{m\in\mathbb{Z}\ |\ f=q^{v}\tau_x^mg, \hbox{ for some }v\in\bZ\}.
%\]
%It is easy to check that unless $f,g\in\bK(x)$ are monomials of the same degree, the $q$-dispersion set $qDisp(f,g)$ is finite and computable. If $qDisp(f,g)$ is not empty, we know $f$ and $g$ are in the same $\la\tau_x\ra$-\emph{orbit}.
%

By Theorem~\ref{th:q-shit equivalent} and partial fraction decomposition w.r.t $y$, any $f\in\bK(x,y)$ can be uniquely decomposed into
\begin{equation}\label{EQ:q-decompose}
f=\mu+yp_1(y)+\frac{p_2(y)}{y^s}
   +\sum_{i=1}^m\sum_{j=1}^{n_i}\sum_{\ell=0}^{k_{i,j}}
          \frac{a_{i,j,\ell}(y)}{\tau_y^{\ell}d_{i}^j(y)},
\end{equation}
where $\mu\in \bK(x)$, $s,m,n_i,k_{i,j}\in \mathbb{N}$, $p_1,p_2,a_{i,j,\ell},d_{i}\in \bK(x)[y]$, $\deg_y(p_2)<s$, $\deg_y(a_{i,j,\ell})<\deg_y(d_i)$, and $d_i\neq y$ are monic irreducible polynomials in distinct $\la\tau_y\ra$-orbits.

For any $n\in\bZ\setminus\{0\}$, as $q^n\neq 1$, we get
\begin{equation}\label{eq:luolang}
y^n=\tau_y\left(\frac{y^n}{q^n-1}\right)-\frac{y^n}{q^n-1}.
\end{equation}
Since $p_1,p_2\in \bK(x)[y]$ with $\deg_y(p_2)<s$, we know both $yp_1$ and $\frac{p_2}{y^s}$ are $\tau_y$-summable in $\bK(x,y)$.
Then the fact that $\tau_y$ preserves the $q$-shift equivalence shows $f$ is $\tau_y$-summable if and only if both $\mu$ and $\sum_{i=1}^m\sum_{j=1}^{n_i}\sum_{\ell=0}^{k_{i,j}}
          \frac{a_{i,j,\ell}(y)}{\tau_y^{\ell}d_{i}^j(y)}$ are $\tau_y$-summable.
Apparently $\mu$ is $\tau_y$-summable only when $\mu=0$.
Thus we only need to consider the $q$-summabiliy of
\begin{equation}\label{form1}
r=\sum_{i=1}^m\sum_{j=1}^{n_i}\sum_{\ell=0}^{k_{i,j}}
\frac{a_{i,j,\ell}}{\tau_y^{\ell}d_{i}^j}.
\end{equation}
\begin{defn}
Let $f\in\bK(x,y)$ be of the form~\eqref{form1}. The sum $\sum_{\ell=0}^{k_{i,j}}\tau_y^{-\ell}a_{i,j,\ell}$
is called the
\emph{$q$-polynomial residue}of $f$ at the $\la\tau_y\ra$-orbit of $d_i$ of multiplicity $j$,
denoted by $qres_y(f,d_i,j)$.
\end{defn}

Note that each summand in Equation~\eqref{form1} is of the form
$a/\tau_y^{\ell} d^j$, which can be transformed by the following Remark.
\begin{rem}\label{reduceorbit}
Suppose $a,d\in\bK(x)[y]$, $j\in\mathbb{N}$ and $\ell\in\mathbb{Z}$. Then we have
\begin{equation}\label{EQ:qunitransformation}
\frac{a}{\tau_y^{\ell}d^j} -\frac{\tau_y^{-\ell}a}{d^j}= \tau_y(g)-g,
\end{equation}
where
\[
g = \begin{cases}
  {
   \sum\limits_{i=0}^{\ell-1}} \frac{\tau_y^{i-\ell}(a)}{\tau_y^i(d^j)},
         & \mbox{if $\ell \ge 0$}, \\[15pt]
  -\sum\limits_{i=0}^{-\ell-1} \frac{\tau_y^i(a)}{\tau_y^{\ell+i}(d^j)},
         & \mbox{if $\ell < 0$
  }.
\end{cases}
\]
This transformation will be used frequently in this paper.
\end{rem}

The following lemma shows if a rational function in $\bK(x,y)$ is $\tau_y$-summable, then we can always rewrite it as the $q$-difference of a rational function in some specific form.

\begin{lem}\label{lemma1}
Suppose a rational function $f\in \bK(x,y)$ of the form \eqref{form1} is $\tau_y$-summable, then there
exists $g\in \bK(x,y)$ such that
$
f=\tau_yg-g
$
and
$g$ is of the form
\begin{equation}\label{eq:reduced form}
g=\sum_{i=1}^{m}\sum_{j=1}^{n_i}\sum_{\ell=0}^{k_{i,j}'}
\frac{c_{i,j,\ell}}{\tau_y^{\ell}d_i^j},
\end{equation}
where $k_{i,j}'\in\bN$, $c_{i,j,\ell},d_i\in\bK(x)[y]$, $d_i$ is monic and $\deg_y(c_{i,j,\ell})<\deg_y(d_i)$.
\end{lem}

\pf Since $f$ is $\tau_y$-summable, there exists $h\in \bK(x,y)$ such that $f=\tau_yh-h$.
Decomposing $h$ according to partial fraction decomposition w.r.t $y$.
Then collecting all fractions whose denominator is of the form $d^j$ with $d\sim_{\la\tau_y\ra}d_i$ into $g$.
Thus we can rewrite $h$ as
$
 h=g+h_1,
$
where $g,h_1\in\bK(x,y)$ and $g$ is of the form
\[
  g=\sum_{i=1}^{m}\sum_{j=1}^{n_i}\sum_{\ell=-v_{i,j}}^{k_{i,j}'}
\frac{c_{i,j,\ell}}{\tau_y^{\ell}d_i^j}
\quad\text{with}\quad v_{i,j},k_{i,j}'\in\bN, c_{i,j,\ell}\in\bK(x)[y].
\]

Taking $h=g+h_1$ into $f=\tau_yh-h$, we deduce that
\[
f-(\tau_yg-g)=\tau_yh_1-h_1.
\]
Considering the denominator of the both sides of the above equation, the choice of $g$ together with the fact that $\tau_x$ preserves the $q$-shift equivalence and multiplicities lead to
\[
f-(\tau_y g-g)=\tau_yh_1-h_1=0.
\]
Next we will prove that all $q$-shifts in the denominator of $g$ are in fact
nonnegative.
Fix $1\leq i\leq m$, $1\leq j\leq n_i$ and let $v_{i,j}=\min\{\ell\in\bZ|c_{i,j,\ell}\neq0\}$.
We claim that $v_{i,j}\geq 0$.
Otherwise $\tau_y^{v_{i,j}}d_i^j$ appears in the denominator of $g$ but not $\tau_yg$ or $f$, which is impossible since $f=\tau_yg-g$.
The arbitrariness of $i,j$ shows $g$ is of the form~\eqref{eq:reduced form}.
\qed

We now ready to present a criterion on the $q$-summability of univariate rational functions via $q$-polynomial residues.

\begin{thm}\label{main1}
Let $f\in \bK(x,y)$ be a rational function of the form \eqref{form1}, then $f$ is $\tau_y$-summable if and only if $qres_y(f,d_i,j)=0$ for any $1\leq i\leq m$ and $1\leq j\leq n_i$.
\end{thm}
\pf For the necessity, if $f$ is $\tau_y$-summable.
Then there exists $g\in \bK(x,y)$ such that $f=\tau_yg-g$ and $g$ is of the form
$
g=\sum_{i=1}^{m}\sum_{j=1}^{n_i}\sum_{\ell=0}^{k_{i,j}^{'}}
\frac{c_{i,j,\ell}}{\tau_y^{\ell}d_i^j}.
$
Then we have
\[
f=\sum_{i=1}^{m}\sum_{j=1}^{n_i}\sum_{\ell=0}^{k_{i,j}^{'}}
                         \frac{\tau_yc_{i,j,\ell}}{\tau_y^{\ell+1}d_i^j}
   -\sum_{i=1}^{m}\sum_{j=1}^{n_i}\sum_{\ell=0}^{k_{i,j}^{'}}
                         \frac{c_{i,j,\ell}}{\tau_y^{\ell}d_i^j}.
\]
Considering the $q$-polynomial residues of both sides of the above equation, we get
\[
qres_y(f,d_i,j)=\sum_{\ell=0}^{k_{i,j}'}\tau_y^{-\ell}c_{i,j,\ell}
                          -\sum_{\ell=0}^{k_{i,j}'}\tau_y^{-\ell}c_{i,j,\ell}=0.
\]
For the sufficiency, if
$qres_y(f,d_i,j)=\sum\limits_{\ell=0}^{k_{i,j}}\tau_y^{-\ell}a_{i,j,\ell}=0$ for all $i,j$, then
\begin{align}\label{decompose}
f & =\sum_{i=1}^m\sum_{j=1}^{n_i}\sum_{\ell=0}^{k_{i,j}}
                 \frac{a_{i,j,\ell}}{\tau_y^{\ell}d_{i}^j}\nonumber\\
  &=\sum_{i=1}^m\sum_{j=1}^{n_i}\sum_{\ell=0}^{k_{i,j}}
    \left( \frac{a_{i,j,\ell}}{\tau_y^{\ell}d_{i}^j}-
             \frac{\tau_y^{-\ell}a_{i,j,\ell}}{d_i^j}\right)
     +\sum_{i=1}^m\sum_{j=1}^{n_i}
         \frac{qres_y(f,d_i,j)}{d_i^j}\nonumber\\
  &=\sum_{i=1}^m\sum_{j=1}^{n_i}\sum_{\ell=0}^{k_{i,j}}
              \left(\tau_yg_{i,j,\ell}-g_{i,j,\ell}\right)\nonumber\\
  &=\tau_yg-g
\end{align}
%%\begin{align*}
%%f & =\sum_{i=1}^m\sum_{j=1}^{n_i}\sum_{\ell=0}^{k_{i,j}}\frac{a_{i,j,\ell}(x)}{d_{i}^j(xq^{\ell})}\\
%%  &=\sum\limits_{i,j,\ell}\left(\tau_x^{\ell}\frac{a_{i,j,\ell}(xq^{-{\ell}})}{d_i^j(x)}-
%%     \frac{a_{i,j,\ell}(xq^{-\ell})}{d_i^j(x)}\right)
%%     +\sum_{i,j}\frac{\sum_{\ell=0}^{k_{i,j}}a_{i,j,\ell}(xq^{-\ell})}{d_i^j(x)}\\
%%  &=\sum\limits_{i,j,\ell}\left(\tau_xg_{i,j,\ell}(x)-g_{i,j,\ell}(x)\right)\\
%%  &=\tau_xg-g
%%\end{align*}
where $g_{i,j,\ell}$ is given according to Remark \ref{reduceorbit} and $g=\sum_{i=1}^m\sum_{j=1}^{n_i}\sum_{\ell=0}^{k_{i,j}}g_{i,j,\ell}$.
This completes the proof.\qed

It should be noticed that, Identity~\eqref{decompose} also shows any rational function $f\in\bK(x,y)$ can be decomposed into $f=\tau_yg-g+r$, where $g\in \bK(x,y)$, $r=\sum\limits_{i,j}\frac{qres_y(f,d_i,j)}{d_i^j}$ with $d_i$ in distinct $\la\tau_y\ra$-orbits.

\medskip\noindent {\bf Univariate qSummability.}
Given a rational function $f\in\bK(x,y)$,
decide whether $f$ is $\tau_y$-summable.
If so, compute a $g\in \bK(x,y)$  such that
\[f = \tau_y(g)-g.\]

\begin{enumerate}
\item Rewrite $f$ into the form \eqref{EQ:q-decompose}, that is
\begin{equation*}
f=\mu+yp_1(y)+\frac{p_2(y)}{y^s}
   +\sum_{i=1}^m\sum_{j=1}^{n_i}\sum_{\ell=0}^{k_{i,j}}
          \frac{a_{i,j,\ell}(y)}{\tau_y^{\ell}d_{i}^j(y)}.
\end{equation*}
\item If $\mu\neq 0$ or $\sum_{\ell=0}^{k_{i,j}}\tau_y^{-\ell}a_{i,j,\ell}\neq 0$ for some $1\leq i\leq m$, $1\leq j\leq n_i$, then return
    ``Not $\tau_y$-summable".	
\item By Identity \eqref{eq:luolang}, obtain a $g$ such that
$yp_1(y)+\frac{p_2(y)}{y^s}=\tau_y(g)-g$.
\item Set $g_{i,j,\ell} =\sum\limits_{k=0}^{\ell-1} \frac{\tau_y^{k-\ell}(a_{i,j,\ell})}{\tau_y^k(d_{i}^j)}$ if $\ell \ge 0$, $g_{i,j,\ell} =-\sum\limits_{k=0}^{-\ell-1} \frac{\tau_y^k(a_{i,j,\ell})}{\tau_y^{\ell+k}(d_i^j)}$  if $\ell < 0$.
\item
Update
$g=g+\sum_{i=1}^m\sum_{j=1}^{n_i}\sum_{\ell=0}^{k_{i,j}}g_{i,j,\ell}$, and then return $g$.
\end{enumerate}
\subsection{Bivariate case}
The goal of this subsection is to decide whether a given rational function $f\in\bK(x,y)$ is $(\tau_x,\tau_y)$-summable.

Viewing $f\in \bK(x, y)$ as a rational function of $y$ over $\bK(x)$, we have the partial fraction decomposition
\begin{equation}\label{EQ:qpfd}
f = \mu(x) +yp_1(x,y)+\frac{p_2(x,y)}{y^s}+\sum_{i,j} \frac{a_{i, j}(x,y)}{d_i^j(x,y)},
\end{equation}
where $s\in\mathbb{Z}$, $\mu(x)\in\bK(x)$, $p_1,p_2\in \bK(x)[y]$ with $\deg_y(p_2)<s$, $a_{i, j}\in \bK(x)[y]$, $d_i\in\bK[x,y]\setminus\{y\}$ are monic and irreducible with  $\deg_y(a_{i,j})<\deg_y(d_i)$.

Since $p_1,p_2\in \bK(x)[y]$ with $\deg_y(p_2)<s$, we know both $yp_1$ and $p_2/y^s$ are $\tau_y$-summable.
Notice that for any $a/\tau_x^m\tau_y^nd^j$ with $m,n\in\mathbb{Z}$, $j\in\mathbb{N}$, $a\in\bK(x)[y]$ and $d\in\bK(x,y)$, applying the transformation~\eqref{EQ:qunitransformation} with respect to $x$ and subsequently
with respect to $y$ yields
\[
  \frac{a}{\tau_x^m\tau_y^nd^j}=\tau_xg-g+\tau_yh-h+\frac{\tau_x^{-m}\tau_y^{-n}a}{d^j}
\]
for some $g,h\in\bK(x,y)$.
Repeating the above transformation, we will get the following decomposition.

\begin{lem}\label{DecomposedForm}
Let $f\in\bK(x,y)$.
Then $f$ can be decomposed into
\[
f=\tau_xg-g+\tau_yh-h+r
\]
with $g,h\in \bK(x,y)$ and $r$ being of the form
\begin{equation}\label{form of r}
r=\mu+\sum_{i=1}^m\sum_{j=1}^{n_i}\frac{a_{i,j}}{d_i^j},
\end{equation}
where $\mu\in\bK(x)$, $a_{i,j}\in \mathbb{K}(x)[y]$, $d_i\in \mathbb{K}[x,y]\setminus\{y\}$ are monic irreducible polynomials, $\deg_y(a_{i,j})<\deg_y(d_i)$,
$d_i \mbox{ and } d_{i'} \mbox{ are not }\la\tau_x,\tau_y\ra\mbox{-equivalent}$
for any $1\leq i, i' \leq m$ and $i\neq i'$.
\end{lem}
Lemma~\ref{DecomposedForm} shows the $(\tau_x,\tau_y)$-summability problem of $f$ is equivalent to that of $r$, and thus is equivalent to determine the $q$-summability of each summand of $r$, as will be proved by the following lemma.

\begin{lem}\label{decompostion to fractions}
Let $r$ be of form~\eqref{form of r}.
Then $r$ is $(\tau_x,\tau_y)$-summable if and only if $\mu$ is $\tau_x$-summable and $a_{i,j}/d_i^j$ is $(\tau_x,\tau_y)$-summable for any $1\leq i\leq m,1\leq j\leq n_i$.
\end{lem}
\pf The sufficiency follows from the linearity of $\tau_x$ and $\tau_y$.
It suffices to prove the necessity.
Suppose $r$ is $(\tau_x,\tau_y)$-summable.
Then it is easy to see both $\mu$ and $\sum_{i=1}^m\sum_{j=1}^{n_i}\frac{a_{i,j}}{d_i^j}$ are $q$-summable.
Then there exist $g,h\in\bK(x,y)$ such that
\[
  \sum_{i=1}^m\sum_{j=1}^{n_i}\frac{a_{i,j}}{d_i^j}=\tau_xg-g+\tau_yh-h.
\]
Decompose $g,h$ as
\[
  g=\sum_{i=1}^{m}\frac{A_i}{B_i}+g_1 \mbox{ and }
  h=\sum_{i=1}^{m}\frac{C_i}{D_i}+h_1,
\]
where $A_i, B_i, C_i, D_i\in \bK[x,y]$,
$\deg_y(A_i)<\deg_y(B_i), \deg_y(C_i)<\deg_y(D_i)$ and $g_1,h_1\in \bK(x,y)$ such that $B_i$ (\emph{resp.} $D_i$) contains exactly all irreducible factors in the denominator of $g$ (\emph{resp.} $h$) which are $\la\tau_x,\tau_y\ra$-equivalent to $d_i$, while the denominators of $g_1$ (resp. $h_1$) contains no such factors.
As $\tau_x,\tau_y$ preserve the $\la\tau_x,\tau_y\ra$-equivalence, we have
\[
  r_i=\sum_{j=1}^{n_i}\frac{a_{i,j}}{d_i^j}
     =\tau_x\frac{A_i}{B_i}-\frac{A_i}{B_i}
      +\tau_y\frac{C_i}{D_i}-\frac{C_i}{D_i}\quad\text{for any }1\leq i\leq m,
\]
and thus $r_i$ is $(\tau_x,\tau_y)$-summable for any $1\leq i\leq m$.

By the observation that $\tau_x$ and $\tau_y$ preserve the multiplicities of irreducible polynomials, we know $r_i$ is $(\tau_x,\tau_y)$-summable if and only if $a_{i,j}/d_i^j$ is $(\tau_x,\tau_y)$-summable for each $j$, which concludes the proof.\qed

\begin{rem}\label{lm:denom}
From the above proof, one can see if $r$ is $(\tau_x,\tau_y)$-summable, we can always find $g,h$ so that $r=\tau_xg-g+\tau_yh-h$ and denominators of $g,h$ contain only irreducible polynomials which are $\la\tau_x,\tau_y\ra$-equivalent to the $d_i$. Such $g,h$ will be refered to as $g,h$ in \emph{reduced form}.
\end{rem}
Since the $q$-summability problem of $\mu\in\bK(x)$ has been solved in the previous subsection, by Lemma~\eqref{decompostion to fractions}, we only need to determine the $q$-summability of fraction $a/d^j$.
\begin{thm}\label{MainTheorem}
Suppose $a/d^j\in \bK(x,y)$, where $j\in \mathbb{N}\setminus\{0\}$, $d\in \bK[x,y]\setminus\{y\}$ is irreducible, $a\in \mathbb{K}(x)[y]\setminus\{0\}$, and
$\deg_y(a)<\deg_y(d)$.
Then $a/d^j$ is $(\tau_x,\tau_y)$-summable if and only if
\begin{itemize}
\item [\rm (1)]
there exist integers $t,\ell,v$ with $t\neq 0$ such that
\begin{equation}\label{condition1}
\tau_x^{t}d=q^{v}\tau_y^{\ell}d,
\end{equation}
\item [\rm (2)]
for the smallest positive integer $t$ satisfying \eqref{condition1}, we have
\begin{equation}\label{condition2}
a=q^{-vj}\tau_x^{t}\tau_y^{-\ell}p-p,
\end{equation}
for some $p\in \mathbb{K}(x)[y]$ with $\deg_y(p)<\deg_y(d)$.
\end{itemize}
\end{thm}
\pf For the sufficiency, since conditions \eqref{condition1} and \eqref{condition2} hold, let $g=\sum\limits_{k=0}^{t-1}\frac{\tau_x^kp}{\tau_x^kd^j}$.
Then
\begin{align}\label{constructed}
  \frac{a}{d^j}-(\tau_xg-g)
=  \frac{a+p}{d^j}
    -\frac{\tau_x^{t}p}{\tau_x^{t}d^j}
=  \tau_y^{-\ell}\left(\frac{\tau_x^tp}{\tau_x^t{d^j}}\right)
   -\frac{\tau_x^tp}{\tau_x^t{d^j}}
=  \tau_yh-h,
\end{align}
where $h$ is given according to Remark \ref{reduceorbit}.
%Then we have found $g,h$, such that
%$\frac{a}{d^j}=\tau_xg-g+\tau_yh-h.$
%Hence $a/d^j$ is $(\tau_x,\tau_y)$-summable.

For the necessity, suppose $a/d^j$ is $(\tau_x,\tau_y)$-summable.
Remark~\ref{lm:denom} shows we can always find $g,h$ in reduced forms such that
\begin{equation}\label{eq:summable}
\frac{a}{d^j}=\tau_xg-g+\tau_yh-h.
\end{equation}
Thus we can assume   $g=\sum\limits_{i=m_0}^m\sum\limits_{k=n_0}^n\frac{b_{i,k}}{\tau_x^i\tau_y^k d^j}$, where $m,m_0,n,n_0\in\bZ$, $b_{i,k}\in\bK(x)[y]$.
By transformation \eqref{EQ:qunitransformation}, we can rewrite $g$ as
\begin{align}\label{eq:g}
g=\sum\limits_{i=m_0}^m\sum\limits_{k=n_0}^n\left\{\tau_y^k\left(\frac{\tau_y^{-k}b_{i,k}}{\tau_x^i d^j}\right)
   -\frac{\tau_y^{-k}b_{i,k}}{\tau_x^i d^j}
   +\frac{\tau_y^{-k}b_{i,k}}{\tau_x^i d^j}\right\}
 =\tau_yg_1-g_1+\sum\limits_{i=m_0}^m\frac{c_i}{\tau_x^id^j},
\end{align}
where $g_1\in\bK(x,y)$, $c_i=\sum\limits_{k=n_0}^n\tau_y^{-k}b_{i,k}$.
Note that at least one of $c_i$ is nonzero,
otherwise \eqref{eq:summable} and Theorem \ref{main1} lead to
\[
  a=qres_y(a/d^j,d,j)=qres_y(\tau_xg-g,d,j)=0,
\]
which contradicts the fact that $a$ is nonzero.

We claim that there exists nonzero integer $t$ such that $\tau_x^td$ is  $\la\tau_y\ra$-equivalent to $d$.
We prove this claim by contradiction. \
Suppose $\tau_x^td$ is not $\la\tau_y\ra$-equivalent to $d$ for any nonzero integer $t$.
If $c_i\neq 0$ for some $i\geq 0$, let $I=\max\{i|c_i\neq 0\}$.
Then by assumption, $\tau_x^{I+1+i}(d)\sim_{\la\tau_y\ra}\tau_x^{I+1}(d)$ if and only if $i=0$.
Computing the $q$-polynomial residues of both sides of \eqref{eq:summable} at the $\la\tau_y\ra$-orbit of $\tau_x^{I+1}d$ of multiplicity $j$, we obtain
\begin{align*}
0=qres_y(a/d^j,\tau_x^{I+1}d,j)
 =qres_y\left(\tau_x\left(\sum\limits_{i}\frac{c_i}{\tau_x^id^j}\right)
              -\sum\limits_{i}\frac{c_i}{\tau_x^id^j},\tau_x^{I+1}d,j\right)
 =\tau_xc_I,
\end{align*}
which is impossible as $c_I\neq 0$.
If $c_i\neq 0$ for some $i<0$, let $i_0=\min\{c_i\neq 0\}$.
Similarly, by computing the $q$-polynomial residue of both sides of \eqref{eq:summable} at the $\la\tau_y\ra$-orbit of $\tau_x^{i_0}d$ of multiplicity $j$, we obtain
\begin{align*}
0=qres_y(a/d^j,\tau_x^{i_0}d,j)
 =qres_y\left(\tau_x\left(\sum\limits_{i}\frac{c_i}{\tau_x^id^j}\right)
              -\sum\limits_{i}\frac{c_i}{\tau_x^id^j},\tau_x^{i_0}d,j\right)
 =-c_{i_0},
\end{align*}
which is also a contradiction.
This completes the proof of the claim.

Suppose $t$ is the smallest positive integer satisfying
$
\tau_x^td=q^{v}\tau_y^{\ell}d.
$
By \eqref{eq:g} and Remark \ref{reduceorbit}, we can rewrite $g$ as
\begin{equation}\label{eq:g-decomposed}
  g=\tau_yg_2-g_2+\sum_{i=0}^{t-1}\frac{\tau_x^ip_i}{\tau_x^id^j},
\end{equation}
where $g_2\in\bK(x,y)$ and $p_i\in\mathbb{K}(x)[y]$.

Substituting \eqref{eq:g-decomposed} into \eqref{eq:summable} and compute the $q$-polynomial residue of both sides of \eqref{eq:summable} at the $\la \tau_y \ra$-orbit of $\tau_x^{i}d$ for $i=0,\ldots,t-1$, we derive
\[
  a=q^{-vj}\tau_x^t\tau_y^{-\ell}p_{t-1}-p_0,\
  0=p_0-p_1,\
  0=p_1-p_2,\
  \ldots,   \
  0=p_{t-2}-p_{t-1}.
\]
Letting $p=p_0=\cdots=p_{t-1}$ leads to
$a=q^{-vj}\tau_x^t\tau_y^{-\ell}p-p$,
which completes the proof.\qed

It is easy to check that the above proof can be adapted to prove Theorem 3.3 in \cite{HouWang2015}, where a much more complicated proof was presented.
\section{Solving the $q$-difference equation}\label{equation}
After all the discussions, we know the $q$-summability problem of a given rational function can be reduced to that of simple fractions $\frac{a(x,y)}{b(x)d^j(x,y)}$, where $j\in\bN\setminus\{0\}$, $a\in\bK[x,y]\setminus\{0\}$, $b\in\bK[x]$, $d\in\bK[x,y]\setminus\{y\}$ is irreducible and $\deg_y(a)<\deg_y(d)=\lambda$.

To determine whether $\frac{a(x,y)}{b(x)d^j(x,y)}$ is $(\tau_x,\tau_y)$-summable, according to Theorem~\ref{MainTheorem}, firstly we need to settle the $q$-shift equivalence testing problem, which has been solved in section \ref{q-dispersion}.
If condition~\eqref{condition1} is not satisfied, then $\frac{a(x,y)}{b(x)d^j(x,y)}$ is not $(\tau_x,\tau_y)$-summable.
Next we always assume condition~\eqref{condition1} holds, and $m$ is the smallest positive integer such that
$\tau_x^m(d) = q^v\tau_y^n(d)$ for some $n,v\in\bZ$.
Then condition~\eqref{condition2} shows we need to solve the following $q$-difference equation
\begin{equation}\label{EQ:qdenom}
\frac{a}{b}=q^{-vj}\tau_x^m\tau_y^{-n}p-p, \quad m,n\in\bZ\text{ and }m>0
\end{equation}
for $p\in\mathbb{K}(x)[y]$ and $\deg_y(p)<\lambda$.
Once we find such a $p$, let
\begin{equation}\label{eq:tranform}
g=\sum_{k=0}^{m-1}\tau_x^k(\frac{p}{d^j})\quad \text{and} \quad
h=\left\{
     \begin{array}{ll}
       -\sum_{k=0}^{n-1}\frac{q^{-vj}\tau_x^m\tau_y^{k-n}p}{\tau_y^kd^j}, & \hbox{when $n\geq 0$;} \\
       \sum_{k=0}^{-n-1}\frac{q^{-vj}\tau_x^m\tau_y^kp}{\tau_y^{k+n}d^j}, & \hbox{when $n<0$.}
     \end{array}
   \right.
\end{equation}

Equation \eqref{constructed} leads to $\frac{a}{bd^j}=\tau_xg-g+\tau_yh-h$.

In this section, we present an algorithm on solving equation \eqref{EQ:qdenom}.
%\begin{equation}\label{EQ:qdenom}
%\frac{a}{b}= c\tau_x^m \tau_y^{-n} p - p,
%\end{equation}
%where $m,n$ are given integers with $m>0$ and $c$ is a nonzero constant.
Noting that $\deg_y(p)<\lambda$, we may assume
\[
p = p_0(x) + p_1(x)y + \cdots + p_{\lambda-1}(x)y^{\lambda-1},\ p_i(x)\in\bK(x).
\]
Comparing the coefficients of like powers of $y$ on both sides of \eqref{EQ:qdenom}, we obtain a system of $q$-difference equations on $p_i(x)$.
%Abramov and coauthors have presented algorithms for finding a universal denominator for the system (see, for example \cite{AbramovPaulePetkovsek1998,Abramov1999rational}).
Next, we will provide an algorithm for finding a universal denominator for the system.

Recall that a $q$-Gosper representation of a rational function $r(x)\in\bK(x)$ is given as $r(x)=\frac{A(x)}{B(x)}\frac{C(xq)}{C(x)},$
where $A(x)$, $B(x)$, $C(x)$ are polynomials, and
$
\gcd(A(x),B(xq^{h}))=1 \text{ for }\forall h\in\bN.
$
Substituting $q^m$ for $q$, we obtain an $m$-fold $q$-Gosper representation $(A(x),B(x),C(x))$ of $r(x)$, that is
\begin{equation}\label{m-fold gosper}
r(x)=\frac{A(x)}{B(x)}\frac{C(xq^m)}{C(x)}\text{ and }
\gcd(A(x),B(xq^{hm}))=1, \forall h\in\bN.
\end{equation}
The following theorem shows that a universal denominator of $p$ can be deduced from an $m$-fold $q$-Gosper representation of $b(x)/b(xq^m)$.

\begin{thm}\label{p_1(x,y)is laurent}
Let $(A(x),B(x),C(x))$ be an $m$-fold $q$-Gosper representation of $\frac{b(x)}{b(xq^m)}$. Then any solution $p(x,y)\in\bK(x)[y]$ to Equation~\eqref{EQ:qdenom} is of the form
\[
p(x,y)=\frac{B(xq^{-m})\hat{p}(x,y)}{b(x)C(x)},
\]
where $\hat{p}(x,y)\in \bK[x,x^{-1},y]$.
\end{thm}
\pf
Rewrite~\eqref{EQ:qdenom} as
\begin{equation}\label{EQ:a=b/b}
a(x,y)=\frac{q^{-vj}b(x)}{b(xq^m)}\tau_x^m\tau_y^{-n}(b(x)p(x,y))-b(x)p(x,y).
\end{equation}
Assume that
\begin{equation}\label{bp}
b(x)p(x,y)=\frac{g(x,y)}{u(x)C(x)},
\end{equation}
where $g(x,y)\in\bK[x,y]$, $u(x)\in\bK[x]$ is a monic polynomial and $\gcd(u(x),g(x,y))=1$.
Substituting~\eqref{bp} into Equation~\eqref{EQ:a=b/b} and noticing that
\begin{equation}\label{b/b m-fold}
\frac{b(x)}{b(xq^m)}=\frac{A(x)}{B(x)}\frac{C(xq^m)}{C(x)},
\end{equation}
we deduce that
\begin{equation}\label{EQ:expan}
a(x,y)B(x)C(x)u(x)u(xq^m)=q^{-vj}A(x)u(x)g(xq^m,yq^{-n})-B(x)u(xq^m)g(x,y).
\end{equation}
Since $\gcd(u(x),g(x,y))=1$, it is readily seen that
\[
u(x)|B(x)u(xq^m) \quad \hbox{ and } \quad u(xq^m)|A(x)u(x).
\]
Next we will show $u(x)$ must be a monomial.
Let $u(x)=x^kh(x)$, for $k\in \bN$ and $h(x)\in \bK[x]$ with $h(0)\neq0$. Then
\begin{equation}\label{eq:h}
h(x)|B(x)h(xq^m) \quad \hbox{ and } \quad h(xq^m)|A(x)h(x).
\end{equation}
If $\gcd(h(x),h(xq^{sm}))=1$ for any $s\in\bZ$, taking $s=0$ yields $h(x)=1$.
If $\gcd(h(x),h(xq^{sm}))\neq 1$ for some $s\in\bZ$, let $\{a_1,a_2,\ldots,a_t\}$ be the set of complex roots of $h(x)=0$, then
\begin{equation}\label{eq:nonprime}
q^{sm}\in \{\frac{a_i}{a_j}:1\leq i,j\leq t\}.
\end{equation}
Since $q^{\ell} \neq 1$ for all
$\ell\in \mathbb{Z}\backslash  \{0\}$, there are only finite $s$ satisfying \eqref{eq:nonprime}.
Let
\[
N=\max\{s\in\bZ:\gcd(h(x),h(xq^{sm}))\neq1\},
\]
and then we have $N\geq0$.
Applying \eqref{eq:h} repeatedly, we obtain
\[h(x)|B(x)B(xq^m)\ldots B(xq^{Nm})
\quad
\text{and}
\quad
h(x)|A(xq^{-m})A(xq^{-2m})\ldots A(xq^{-(N+1)m}).\]
Then \eqref{m-fold gosper} yields $h(x)=1$.
At this stage, we have shown $u(x)$ is a monomial.
Substituting $u(x)=x^k$ into \eqref{EQ:expan} we obtain
\[
a(x,y)B(x)C(x)u(xq^m)=q^{-vj}A(x)g(xq^m,yq^{-n})-B(x)g(x,y)q^{mk},
\]
which leads to
$
B(x)|g(xq^m,yq^{-n}).
$
Assuming $g(x,y)=B(xq^{-m})\bar{p}(x,y)$. Then Identity \eqref{bp} shows
$
p(x,y)=\frac{B(xq^{-m})\bar{p}(x,y)}{x^kb(x)C(x)}.
$
Setting $\hat{p}(x,y)=\frac{\bar{p}(x,y)}{x^k}$ concludes the proof.\qed

Substituting $p(x,y)=\frac{B(xq^{-m})\hat{p}(x,y)}{b(x)C(x)}$ into~\eqref{EQ:a=b/b}, we obtain
\begin{equation}\label{qequation2}
a(x,y)C(x)=q^{-vj}A(x)\hat{p}(xq^m,yq^{-n})-B(xq^{-m})\hat{p}(x,y).
\end{equation}
Notice that $\deg_y(\hat{p})=\deg_y(p)<\lambda$ and $\hat{p}$ is a Laurent polynomial in $x$.
In order to solve~\eqref{qequation2} for $\hat{p}(x,y)$, it suffices to find
bounds on the highest and lowest degrees of $\hat p$ w.r.t. $x$.
Then we can make an ansatz with respect to the coefficients of $x,y$ and obtain a system of linear equations, which can be solved by the classical Gaussian elimination and many other algorithms.
%%
%%knowing the lowest degree can reduce the problem to solving a linear system of $q$-difference equations on the coefficients of the numerator of $\hat{p}$ w.r.t. y.
%%Then a upper bound of $\hat p$ w.r.t. $x$ enables us to solve the obtained linear system by ansatzes.

According to the $m$-fold $q$-Gosper representation~\eqref{b/b m-fold}, we can write
\begin{align*}
& b(x)=\sum_{k=\ell_0'}^{\ell_0}\beta_kx^k,\quad
  A(x)=\sum_{k=\ell_1'}^{\ell_1}a_kx^k, \quad
  B(x)=\sum_{k=\ell_1'}^{\ell_1}b_kx^k,\\
& C(x)=\sum_{k=\ell_2'}^{\ell_2}c_kx^k,\quad
  a(x,y)=\sum_{k=\ell_3'}^{\ell_3}\alpha_k(y)x^k,
  \quad \hat{p}(x,y)=\sum_{k=\ell'}^{\ell}p_k(y)x^k,
\end{align*}
where $\ell,\ell',\ell_i,\ell_i'\in\bZ$ are the highest and lowest degrees
of those polynomials.
Note that
$b_{\ell_1}=a_{\ell_1}q^{m(\ell_0+\ell_2)}$ and $b_{\ell_1'}=a_{\ell_1'}q^{m(\ell_0'+\ell_2')}$
by comparing coefficients of Equation~\eqref{b/b m-fold}.

%%The following theorem provides the bounds on $\deg_x(\hat{p})$.
\begin{thm}\label{upper lower bound}
Suppose $(A(x), B(x), C(x))$ is an $m$-fold $q$-Gosper representation of~$\frac{b(x)}{b(xq^m)}$, $\deg_y(a)<\lambda$ and the above expressions hold. If $\hat{p}(x,y)\in\bK(x)[y]$ is a Laurent polynomial in $x$ which satisfies~\eqref{qequation2}, then we have
\[
\ell\leq \max \left\{\ell_2+\ell_3-\ell_1, \ell_0-\ell_1+\ell_2+\frac{vj}{m}, \ell_0-\ell_1+\ell_2+\frac{vj+n(\lambda-1)}{m}\right\}
\]
and
\[
\ell'\geq \min \left\{\ell_2'+\ell_3'-\ell_1', \ell_0'-\ell_1'+\ell_2'+\frac{vj}{m},
                            \ell_0'-\ell_1'+\ell_2'+\frac{vj+n(\lambda-1)}{m}.\right\}
\]
\end{thm}
\pf Taking \eqref{qequation2} as an equation in $x$ over $\bK[y]$, firstly we try to find an upper bound on $\deg_x(\hat{p})$.
Suppose the leading coefficient of the right hand of \eqref{qequation2} is not canceled, then we arrive at
\begin{equation}\label{EQ:Ncancel}
\ell=\ell_2+\ell_3-\ell_1.
\end{equation}
Otherwise we have
$
q^{-vj}a_{\ell_1}p_{\ell}(yq^{-n})q^{m\ell}-b_{\ell_1}q^{-m\ell_1}p_{\ell}(y)=0.
$
Since $b_{\ell_1}=a_{\ell_1}q^{m(\ell_0+\ell_2)}$, we obtain
\begin{equation}\label{EQ:uppercoe1}
p_{\ell}(yq^{-n})=q^{m(\ell_0-\ell_1+\ell_2-\ell)+vj}p_{\ell}(y).
\end{equation}
If $n=0$, we derive that
\begin{equation}\label{EQ:cancel0}
\ell=\ell_0-\ell_1+\ell_2+\frac{vj}{m}.
\end{equation}
Otherwise \eqref{EQ:uppercoe1} implies that
$
p_{\ell}(y)=\epsilon y^k, \hbox{ for some }\epsilon\in \bK\hbox{ and }0\leq k\leq \lambda-1.
$
Substituting $p_{\ell}(y)=\epsilon y^k$ into \eqref{EQ:uppercoe1}, we get
\begin{equation}\label{EQ:cancel1}
\ell =\ell_0-\ell_1+\ell_2+\frac{vj}{m}+\frac{nk}{m}.
\end{equation}
%%\begin{align*}
%%\ell & =\ell_0-\ell_1+\ell_2-\frac{\log_qc}{m}+\frac{nk}{m}\\
%%    & \leq \max \{\ell_0-\ell_1+\ell_2-\frac{\log_qc}{m}, \ell_0-\ell_1+\ell_2-\frac{\log_qc}{m}+\frac{n(d-1)}{m}\}.
%%\end{align*}
Since $n\in\bZ\setminus\{0\}$ and $0\leq k\leq \lambda-1$, Equation~\eqref{EQ:cancel1} together with \eqref{EQ:Ncancel} and \eqref{EQ:cancel0} lead to
\[
\ell\leq \max \left\{\ell_2+\ell_3-\ell_1, \ell_0-\ell_1+\ell_2+\frac{vj}{m}, \ell_0-\ell_1+\ell_2+\frac{vj+n(\lambda-1)}{m}\right\}.
\]

Considering the coefficient of the lowest degree of the right hand of \eqref{qequation2}, similar discussions conclude the proof.\qed

\noindent {\bf Solving $q$-Difference Equation}
Given $(a/b,m,n,v,j,\lambda)$, where $a\in\bK[x,y]\setminus\{0\}$, $b\in\bK[x]$, $n,v\in\bZ$, $m,j,\lambda\in\bN\setminus\{0\}$ and $\deg_y(a)<\lambda$.
Decide whether $q$-difference equation
$\frac{a}{b}=q^{-vj}\tau_x^m\tau_y^{-n}p-p
$
has a solution $p\in\mathbb{K}(x)[y]$ and $\deg_y(p)<\lambda$.
If so, compute such a $p$.
\begin{enumerate}
\item Compute an $m$-fold $q$-Gosper representation $(A(x),B(x),C(x))$ of $b(x)/b(xq^m)$.
\item Let $\ell$ and $\ell'$ be as given by Theorem~\ref{upper lower bound} and $\hat{p}(x,y)=\sum_{k=\ell'}^{\ell}p_k(y)x^k$.
    Compare the coefficients of $x,y$ both sides of \eqref{qequation2} and solve the obtained linear equations for $p_k$.\\ [1ex]
    \hphantom{Let} If no solution found, return ``No solution". \\[1ex]
    \hphantom{Let} Otherwise, let $\hat{p}$ be one of the solutions.
\item Return ``p=$\frac{B(xq^{-m})\hat{p}(x,y)}{b(x)C(x)}$".
\end{enumerate}

With everything in place, we are now ready to give an algorithm for deciding the $q$-summability.

\noindent {\bf Bivariate qSummability.}
Given a rational function $f\in\bK(x,y)$,
decide whether $f$ is $q$-summable.
If so, compute $g, h\in \bK(x,y)$  such that
\[f = \tau_x(g)-g+\tau_y(h)-h.\]

\begin{enumerate}
\item By partial fraction decomposition and transformation~\eqref{EQ:qunitransformation}, rewrite $f$ as
     \[
       f=\tau_xg-g+\tau_yh-h+r
     \]
     with $r$ of the form~\eqref{form of r}, that is
     $r=\mu+\sum_{i=1}^m\sum_{j=1}^{n_i}\frac{a_{i,j}}{d_i^j}$.
\item Apply algorithm {\bf Univariate qSummability} to $\mu$.\\[1ex]
      \hphantom{Apply} If $\mu$ is not $\tau_x$-summable, then return``Not $(\tau_x,\tau_y)$-summable."\\[1ex]
      \hphantom{Apply} Otherwise, update $g$ to be $g+${\bf Univariate qSummability}$(\mu)$.
\item For $i = 1, \dots, m$ do\\[1ex]
  \hphantom{For}If there do not exist integers $k_i,\ell_i,v_i$ with $t_i\neq 0$ such that
  $\tau_x^{t_i}d_i=q^{v_i}\tau_y^{\ell{_i}}d_i$,
     then return ``Not $(\tau_x,\tau_y)$-summable."\\[1ex]
  \hphantom{For}Otherwise let $m_i$ be the smallest positive integer such that $\tau_x^{m_i}d_i=q^{v_i}\tau_y^{n{_i}}d_i$ and $\lambda_i=\deg_y(d_i)$.

     For $j=1,\dots, n_i$ do\\[1ex]
    \hphantom{For} Applying algorithm {\bf Solving $q$-Difference Equation} to $(a_{i,j},m_i,n_i,v_i,j,\lambda_i)$.\\[1ex]
        \hphantom{apply} If no solution found, then return ``Not $(\tau_x,\tau_y)$-summable."\\[1ex]
        \hphantom{apply} Otherwise, let $p$={\bf Solving $q$-Difference Equation}$(a_{i,j},m_i,n_i,v_i,j,\lambda_i)$. \\[1ex]
        \hphantom{apply} Update  $g=g+\sum_{k=0}^{m_i-1}\tau_x^k(\frac{p}{d_{i}^j})$,
        $h=h-\sum_{k=0}^{n_i-1}\frac{q^{-{v_i}j}\tau_x^{m_i}\tau_y^{k-{n_i}}p}
       {\tau_y^kd_{i}^j}$ when $n_i\geq 0$, 
        \hphantom{apply} $h=h+\sum_{k=0}^{-{n_i}-1}\frac{q^{-{v_i}j}\tau_x^{m_i}
       \tau_y^kp}{\tau_y^{k+n_i}d_{i}^j}$ when $n_i<0$
\item Return $g,h$.
\end{enumerate}
Next we will illustrate how to use our criteria to determine the $q$-summability.
\begin{exam}
Suppose $f$ admits the following partial fraction decomposition
\[
f(x,y)=\frac{y}{xq^2(x^nq^{n}+y^nq^{-n})}+\frac{1}{(x^n+y^n)}+\frac{1}{(x+y)(x+1)}, \]
where $n>1$ is an integer.
Note that $x^nq^{n}+y^nq^{-n}=\tau_x\tau_y^{-1}(x^n+y^n)$. Applying Remark \ref{reduceorbit} repeatedly, we obtain
\[
  f=\tau_x(g_0)-g_0+\tau_y(h_0)-h_0+\frac{x+y}{x(x^n+y^n)}+\frac{1}{(x+y)(x+1)},
\]
where
$g_0=-h_0$
and
$h_0=\frac{-yq^{-1}}{x(x^n+y^nq^{-n})}.$

Thus the $q$-summability of $f$ is equivalent to that of  $r=\frac{x+y}{x(x^n+y^n)}+\frac{1}{(x+y)(x+1)}$.
As $x^n+y^n$ and $x+y$ are not $q$-shift equivalent, Lemma \ref{decompostion to fractions} shows $r$ is $q$-summable if and only if $\frac{x+y}{x(x^n+y^n)}$ and $\frac{1}{(x+y)(x+1)}$ are both $q$-summable.
For $\frac{x+y}{x(x^n+y^n)}$,
it is easy to check that
$
\tau_x(x^n+y^n)=q^n\tau_y^{-1}(x^n+y^n)
$
and $p_n=\frac{q^n}{1-q^n}\cdot (\frac{y}{x}+1)$ is a solution of
$
\frac{x+y}{x}=q^{-n}\tau_x\tau_yp-p.
$
Then identitity~\eqref{eq:tranform} leads to
\begin{equation*}
\frac{x+y}{x(x^n+y^n)}=\tau_x\left(\frac{p_n}{x^n+y^n}\right)-\frac{p_n}{x^n+y^n}
                  +\tau_y\left(\frac{\tau_xp_n}{x^nq^n+y^n}\right)
                  -\frac{\tau_xp_n}{x^nq^n+y^n}.
\end{equation*}
Since $\tau_x(x+y)=q\tau_y^{-1}(x+y)$, we know from Theorem~\ref{MainTheorem} that $\frac{1}{(x+y)(x+1)}$ is  $(\tau_x,\tau_y)$-summable if and only if
\begin{equation}\label{ex:not}
\frac{1}{x+1}=q^{-1}\tau_x\tau_yp(x)-p(x)
\end{equation}
has a solution in $\bK(x)$.
However Theorem~\ref{p_1(x,y)is laurent} together with the fact that $(x+1)/(xq+1)$ admits an $m$-fold Gosper representation $(x+1,xq+1,1)$ imply that $p(x)$ satisfying Equation $\eqref{ex:not}$ must be a Laurent polynomial in $x$, which is impossible since $x+1$ and $x$ are not $q$-shift equivalent.
Thus $r$ hence $f$ is not $(\tau_x,\tau_y)$-summable.
\end{exam}

\noindent \textbf{Acknowledgments.} I would like to thank Shaoshi Chen, Hui Huang and Lixin Du for their constructive comments.
This work was supported by the Natural Science Foundation of Tianjin(19JCQNJC14500) and the National Science Foundation(11871067) of China.

\bibliographystyle{plain}

\bibliography{sum}

\end{document}